\documentclass{article}

\usepackage[english,portuges]{babel}

\usepackage{amsmath,amsfonts,amssymb,amsthm}
\usepackage{graphicx}
\usepackage{url}

\begin{document}

\selectlanguage{english}

\title{Dynamics of Dengue epidemics\\
using optimal control\thanks{Part of the first author's
PhD project carried out at the Universities of Aveiro and Minho.
Submitted to \emph{Mathematical and Computer Modelling} 25/Oct/2009;
accepted for publication, after revision, 22/June/2010.}}

\author{Helena Sofia Rodrigues$^1$\\
\url{sofiarodrigues@esce.ipvc.pt}
\and
M. Teresa T. Monteiro$^2$\\
\url{tm@dps.uminho.pt}
\and
Delfim F. M. Torres$^3$\\
\url{delfim@ua.pt}}

\date{$^1$School of Business Studies\\
  Viana do Castelo Polytechnic Institute\\
  4900-347 Viana do Castelo, Portugal\\[0.3cm]
$^2$Department of Production and Systems\\
  University of Minho\\
  4710-057 Braga, Portugal\\[0.3cm]
$^3$Department of Mathematics\\
  University of Aveiro\\
  3810-193 Aveiro, Portugal}

\maketitle


\begin{abstract}
We present an application of optimal control theory to Dengue epidemics.
This epidemiologic disease is an important theme in tropical countries
due to the growing number of infected individuals. The dynamic model
is described by a set of nonlinear ordinary differential equations,
that depend on the dynamic of the Dengue mosquito, the number of infected individuals,
and the people's motivation to combat the mosquito.
The cost functional depends not only on the costs of medical treatment
of the infected people but also on the costs related to educational and sanitary campaigns.
Two approaches to solve the problem are considered: one using optimal control theory,
another one by discretizing first the problem
and then solving it with nonlinear programming.
The results obtained with \textsf{OC-ODE}
and \textsf{IPOPT} solvers are given and discussed.
We observe that with current computational tools it is easy to obtain,
in an efficient way, better solutions to Dengue problems,
leading to a decrease of infected mosquitoes and individuals
in less time and with lower costs.

\medskip

\textbf{Keywords:} optimal control, Dengue,
nonlinear programming.

\smallskip

\textbf{2010 Mathematics Subject Classification:}
49M25, 49M37, 90C30.

\end{abstract}


\section{Introduction}

Dengue is a mosquito-borne infection, transmitted
by the \emph{Aedes aegytpi} mosquito, that causes
a severe flu-like illness, and sometimes a potentially
lethal complication called Dengue haemorrhagic fever.
Dengue is found in tropical and sub-tropical climates worldwide,
mostly in urban and semi-urban areas. According to the
World Health Organization \cite{WHO}, the incidence of Dengue
has grown dramatically in recent decades. About
40\% of world's population are now at risk.

The aim of the paper is to present a mathematical model
to study the dynamic of the Dengue epidemics,
in order to minimize the investments in disease's control,
since the financial resources are always scarce.
Quantitative methods are applied to the optimization
of investments in the control of the epidemiologic disease,
in order to obtain a maximum of benefits
from a fixed amount of financial resources.
The used model depends on the dynamic of the growing of the mosquito,
but also on the efforts of the public management to motivate
the population to break the reproduction cycle of the mosquitoes
by avoiding the accumulation of still water in open-air recipients
and spraying potential zones of reproduction.

The paper is organized as follows. Section~\ref{sec:2} presents
the dynamic model for Dengue epidemics, where the variables, parameters,
and ordinary differential equations describing the control system, are defined.
In Section~\ref{sec:3} the numerical implementation and the strategies
used to solve the problem are shown. Section~\ref{sec:4} reports
the obtained numerical results. The main conclusions
are then given in Section~\ref{sec:conc}.


\section{Dynamic model}
\label{sec:2}

The Dengue epidemic model described in this paper is based
on the one proposed in \cite{Caetano2001}.
It consists in minimizing
\begin{equation}
\label{cost}
J\left[u_1(\cdot),u_2(\cdot)\right]
=\int_{0}^{t_f}\{\gamma_D x_{3}^{2}(t)
+\gamma_F u_{1}^{2}(t)+\gamma_E u_{2}^{2}(t)\}dt
\end{equation}
subject to the following four nonlinear time-varying state equations
\cite{Caetano2001}:
\begin{equation}
\dot{x}_{1}(t)
=\left[\alpha_{R}\left(1-\mu \sin(\omega t+\varphi)\right)
-\alpha_M - x_4(t)\right]x_1(t)-u_1(t), \label{x1}
\end{equation}
\begin{equation}
\dot{x}_{2}(t)=\left[\alpha_{R}\left(1-\mu \sin(\omega t+\varphi)\right)
-\alpha_M - x_4(t)\right]x_2(t)+\beta \left[ x_1(t)-x_2(t)\right]x_3(t)-u_1(t), \label{x2}
\end{equation}
\begin{equation}
\dot{x}_{3}(t)= -\eta x_3(t)+\rho x_2(t)\left[P-x_3(t)\right], \label{x3}
\end{equation}
\begin{equation}
\dot{x}_{4}(t)= -\tau x_4(t)+\theta x_3(t)+u_2(t), \label{x4}
\end{equation}
where $\dot{x}_i(t) = \frac{dx_{i}(t)}{dt}$, $i = 1, \ldots, 4$.
The notation used in the mathematical formulation
(\ref{cost})-(\ref{x4}) is as follows.

\medskip

\begin{tabular}{lp{10cm}}
\multicolumn{2}{l}{\textbf{State Variables:}}\\
$x_{1}(t)$ & density of mosquitoes \\
$x_{2}(t)$ & density of mosquitoes carrying the virus \\
$x_{3}(t)$ & number of individuals with the disease \\
$x_{4}(t)$ & level of popular motivation to combat mosquitoes (goodwill) \\
\end{tabular}

\medskip

\begin{tabular}{lp{10cm}}
\multicolumn{2}{l}{\textbf{Control Variables:}}\\
$u_{1}(t)$ & investments in insecticides \\
$u_{2}(t)$ & investments in educational campaigns \\
\end{tabular}

\medskip

\begin{tabular}{lp{12.5cm}}
\multicolumn{2}{l}{\textbf{Parameters:}}\\
$\alpha_{R}$ & average reproduction rate of mosquitoes \\
$\alpha_{M}$ & mortality rate of mosquitoes \\
$\beta$ & probability of contact between non-carrier mosquitoes and infected individuals\\
$\eta$ & rate of treatment of infected individuals \\
$\mu$ & amplitude of seasonal oscillation in the reproduction rate of mosquitoes \\
$\rho$ & probability of individuals becoming infected \\
$\theta$ & fear factor, reflecting the increase in the population
willingness to take actions to combat the mosquitoes as a consequence
of the high prevalence of the disease in the specific social environment \\
$\tau$ & forgetting rate for goodwill of the target population \\
$\varphi$ & phase angle to adjust the peak season for mosquitoes \\
$\omega$ & angular frequency of the mosquitoes proliferation cycle,
corresponding to a 52 weeks period\\
$P$ & population in the risk area (usually normalized to yield $P=1$) \\
$\gamma_D$ & the instantaneous costs due to the existence of infected individuals\\
$\gamma_F$ & the costs of each operation of spraying insecticides\\
$\gamma_E$ & the cost associated to the instructive campaigns \\
\end{tabular}

\bigskip

Equation (\ref{x1}) represents the variation of the density of mosquitoes per unit time
to the natural cycle of reproduction and mortality ($\alpha_R$ and $\alpha_M$),
due to seasonal effects $\mu \sin(\omega t+\varphi)$ and to human interference $- x_4(t)$ and $u_1(t)$.
Equation (\ref{x2}) expresses the variation of the density of mosquitoes carrying the virus $x_2$.
The term $\left[\alpha_{R}\left(1-\mu \sin(\omega t+\varphi)\right)-\alpha_M - x_4(t)\right]x_2(t)$
represents the rate of the infected mosquitoes and $\beta \left[ x_1(t)-x_2(t)\right]x_3(t)$
represents the increase rate of the infected mosquitoes due to the possible contact between
the non infected mosquitoes $x_1(t)-x_2(t)$ and individuals with disease denoted by $x_3(t)$.
The dynamics of the infectious transmission is presented in equation (\ref{x3}).
The term $-\eta x_3(t)$ represents the rate of cure and $\rho x_2(t)\left[P-x_3(t)\right]$
represents the rate at which new cases spring up. The factor $\left[P-x_3(t)\right]$
is the number of individuals in the area, that are not infected.
Equation (\ref{x4}) is a model for the level of popular motivation (or goodwill)
to combat the reproductive cycle of mosquitoes.
Along the time, the level of people's motivation changes.
As a consequence, it is necessary to invest in educational
campaigns designed to increase consciousness of the population under risk.
The expression $-\tau x_4(t)$ represents the decay of the people's motivation with time,
due to forgetting. The expression $\theta x_3(t)$ represents the natural
sensibilities of the public due to increase in the prevalence of the disease.

The goal is to minimize the cost functional (\ref{cost}).
This functional includes the social costs related to the existence of ill individuals,
$\gamma_D x_{3}^{2}(t)$, the recourses needed for spraying of insecticides operations,
$\gamma_F u_{1}^{2}(t)$, and for educational campaigns, $\gamma_E u_{2}^{2}(t)$.
The model for the social cost is based on the concept
of goodwill explored by Nerlove and Arrow \cite{Nerlove1962}.

Due to computational issues, the optimal control problem (\ref{cost})-(\ref{x4}),
that is written in the Lagrange form, was converted into an equivalent Mayer problem.
Hence, using a standard procedure (\textrm{cf.}, \textrm{e.g.},
\cite{Lewis1995}) to rewrite the cost functional,
the state vector was augmented by an extra component $x_5$,
\begin{equation}
\label{newx5}
\dot{x}_5(t)=\gamma_D x_{3}^{2}(t)
+\gamma_F u_{1}^{2}(t)+\gamma_E u_{2}^{2}(t),
\end{equation}
leading to the following equivalent terminal cost problem: to minimize
\begin{equation}
\label{x5tf}
I[x_5(\cdot)]=x_5(t_f),
\end{equation}
with given $t_f$, subject to the control system (\ref{x1})-(\ref{x4}) and (\ref{newx5}).


\section{Numerical implementation}
\label{sec:3}

The simulations were carried out using the following normalized numerical values:
$\alpha_{R}=0.20$, $\alpha_{M}=0.18$, $\beta=0.3$,
$\eta=0.15$, $\mu=0.1$, $\rho=0.1$, $\theta=0.05$, $\tau=0.1$, $\varphi=0$,
$\omega=2\pi/52$, $P=1.0$, $\gamma_D=1.0$, $\gamma_F=0.4$, $\gamma_E=0.8$,
$x_1(0)=1.0$, $x_2(0)=0.12$, $x_3(0)=0.004$, and $x_4(0)=0.05$.
These values are available in the paper \cite{Caetano2001}
and were adopted here in order to be able to compare the obtained results
with those of \cite{Caetano2001}. It was considered as final time $t_f=52$ weeks.

Two different implementations were considered.
In the first one, a direct approach was followed through
the use of the \textsf{OC-ODE} optimal control solver.
In the second case, with the employment of a nonlinear solver in mind,
it was necessary to discretize the optimal control problem into a standard
nonlinear programming problem.

The \textsf{OC-ODE} \cite{Matthias2009}, \emph{Optimal Control of Ordinary-Differential Equations},
is a collection of \textsf{Fortran 77} routines for optimal control problems subject
to ordinary differential equations. It uses an automatic direct discretization method
and includes procedures for numerical adjoint estimation and sensitivity analysis.
In our case the formulation used is optimal control problem
(\ref{x1})-(\ref{x4}) and (\ref{newx5})-(\ref{x5tf}).

The \textsf{IPOPT} \cite{Ipopt}, \emph{Interior Point OPTimizer},
is a software package for large-scale nonlinear optimization.
It is written in \textsf{Fortran} and \textsf{C}. \textsf{IPOPT} implements
a primal-dual interior point method and uses a line search strategy based on filter method.
\textsf{IPOPT} can be used from various modeling environments.
In this work, the problem was coded in \textsf{AMPL} \cite{AMPL}
and interfaced to \textsf{IPOPT}.

For the nonlinear strategy, we had to discretize the problem ourselves.
It was selected a first order method: the Euler's scheme \cite{Betts,Edite}.
Higher order discretization schemes can also be considered,
but bringing no advantage \cite{Nosso1}.

It is assumed that the time $t=nh$ moves ahead in uniform steps of length $h$.
So, if a differential equation is written like $\displaystyle\frac{dx}{dt}=f(t,x)$,
it is possible to make a convenient approximation of this:

\begin{center}
\begin{tabular}{l}
$\displaystyle\frac{dx}{dt}=f(t,x)$\\
\\
$\displaystyle\frac{x(t_{n+1})-x(t_n)}{h}\simeq f(t,x)$\\
\\
$\displaystyle
\Rightarrow x_{n+1}\simeq x_n+hf(t_n,x_n)$.
\end{tabular}
\end{center}

This approximation $x_{n+1}$ of $x(t)$ at the point $t_{n+1}$ has an error of order $h^2$.
So it is important to use small enough steps, in order to obtain an accurate solution, in spite of allowing
a computational efficient scheme. The value $h=1/4$ was found to be a good compromise
between precision and efficiency, and was adopted here.
Thus, the optimal control problem was discretized into
the following nonlinear programming problem:
\begin{center}
\begin{tabular}{lll}
minimize & \multicolumn{2}{l}{$x_5(N)$}\\
s.t. &  $x_{1}(i+1)=$ & $x_{1}(i)+h\left\{\left[\alpha_{R}\left(1-\mu
\,\sin(\omega i+\varphi)\right)\right.\right.$\\
& & $\left.\left.-\alpha_M - x_4(i)\right]x_1(i)-u_1(i)\right\}$\\
& $x_{2}(i+1)=$ & $x_{2}(i)+h\left\{\left[\alpha_{R}\left(1-\mu
\,\sin(\omega i+\varphi)\right)-\alpha_M - x_4(i)\right]x_2(i)\right.$\\
& & $\left. + \beta \left[ x_1(i)-x_2(i)\right]x_3(i)-u_1(i)\right\}$\\
&  $x_{3}(i+1)=$ & $x_{3}(i)+h\{-\eta x_3(t)+\rho x_2(t)\left[P-x_3(t)\right]\}$\\
&  $x_{4}(i+1)=$ & $x_{4}(i)+h\{-\tau x_4(t)+\theta x_3(t)+u_2(t)\}$\\
& $x_{5}(i+1)=$ & $x_{5}(i)+h\{\gamma_D x_{3}^{2}(t)+\gamma_F u_{1}^{2}(t)+\gamma_E u_{2}^{2}(t)\}$,\\
& & \\
\multicolumn{3}{l}{where $i \in \{0,\ldots,N-1\}$.}
\end{tabular}
\end{center}

The error tolerance value was $10^{-8}$ using the \textsf{IPOPT} solver.
The \textsf{NEOS} Server \cite{NEOS} platform was used as interface with the solver.
\textsf{NEOS} (\emph{Network Enabled Optimization System}) is an optimization
service that is available through the Internet. There is a large set
of software packages, considered as the state of the art in optimization.
The discretized problem, after a presolve done by the software,
is transformed into 1455 variables, 1243 of which are nonlinear;
and 1039 constraints, 828 of which are nonlinear.

Next section presents the main obtained results.


\section{Computational results}
\label{sec:4}

\begin{figure}[ptbh]
\centering
\begin{minipage}[t]{0.50\linewidth}
\centering
\includegraphics[scale=0.36]{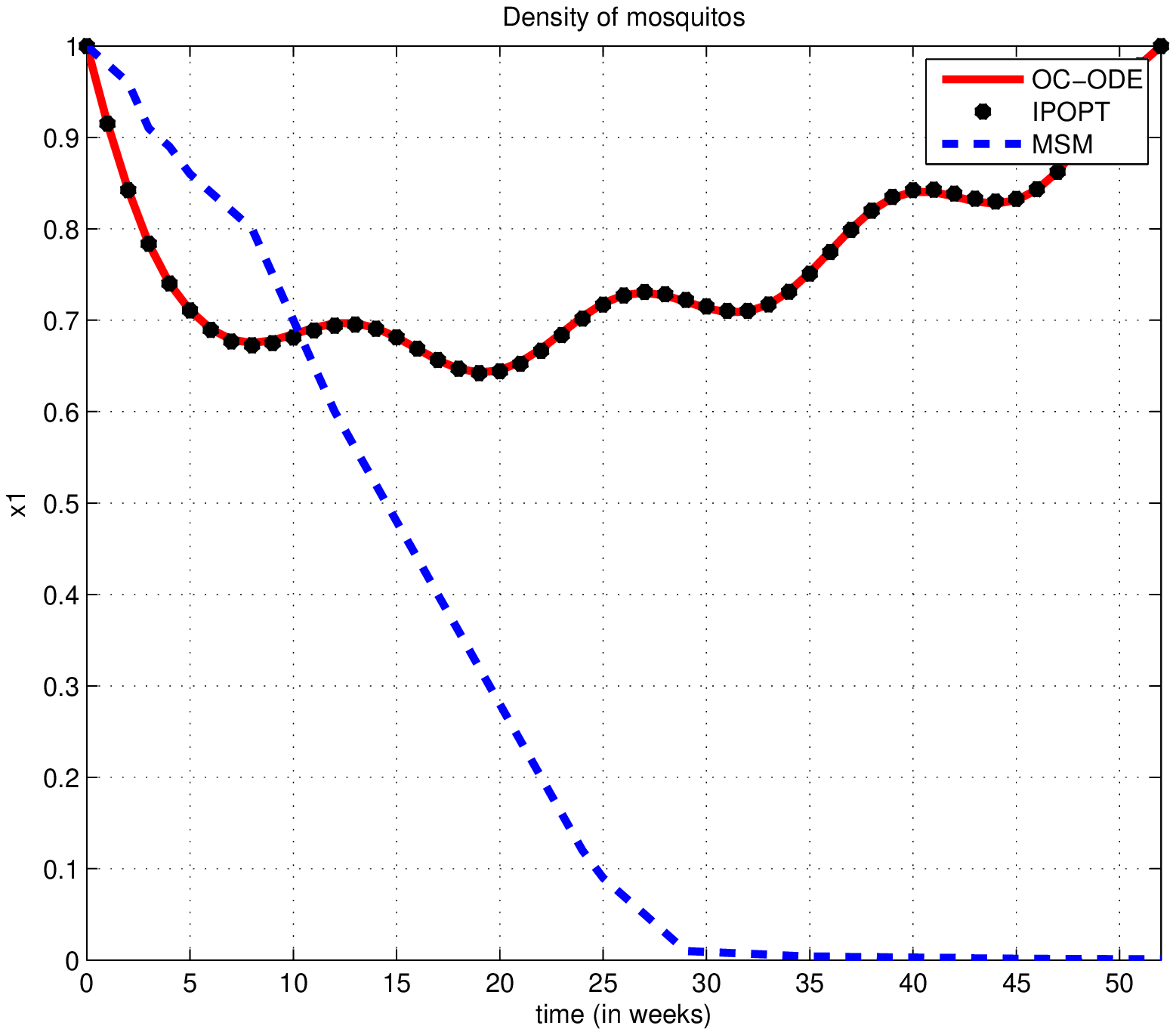}
{\caption{\label{plotx1}  Density of mosquitoes}}
\end{minipage}\hspace*{\fill}
\begin{minipage}[t]{0.50\linewidth}
\centering
\includegraphics[scale=0.36]{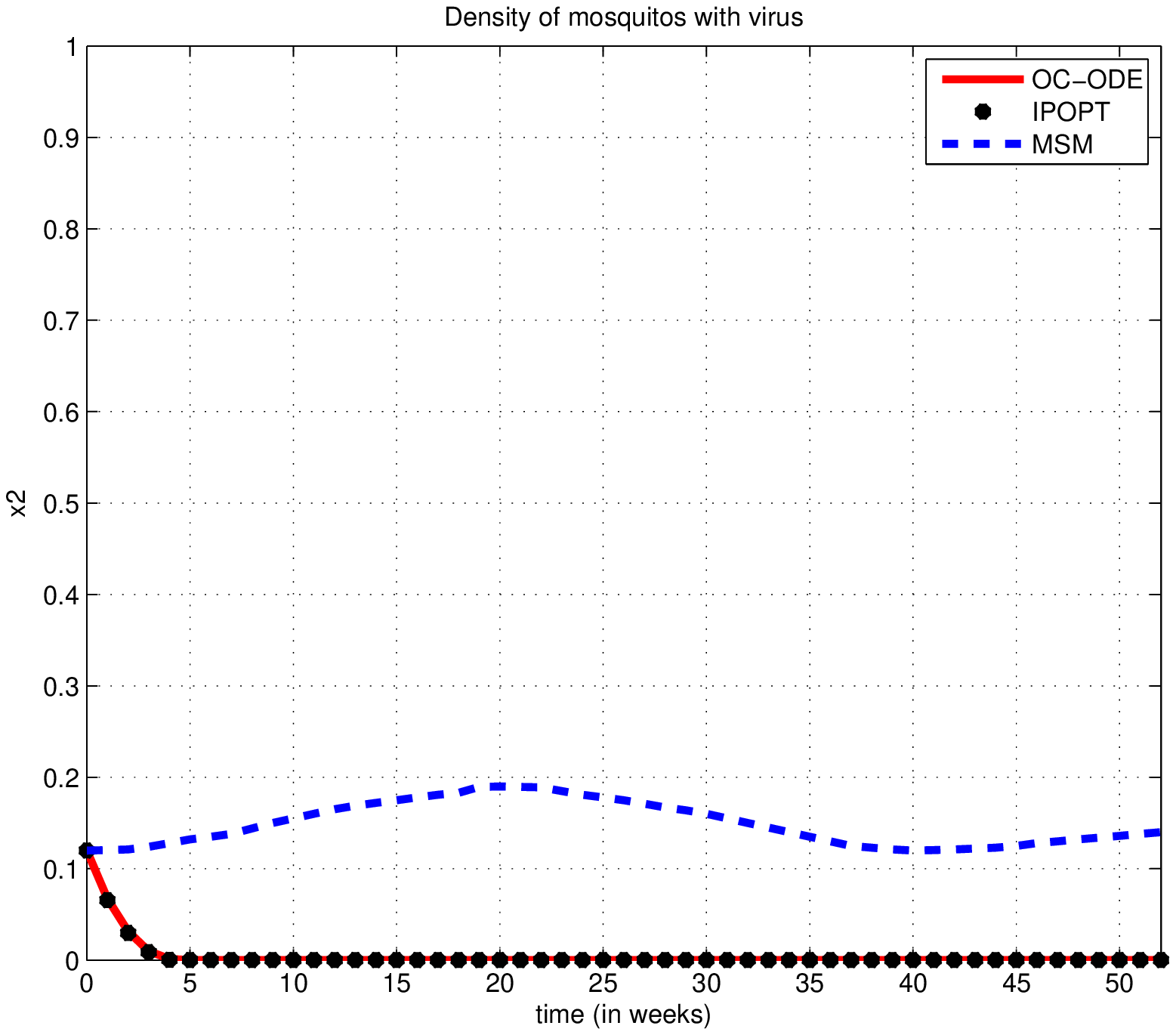}
{\caption{\label{plotx2} \small Infected mosquitoes}} \label{plot}
\end{minipage}
\end{figure}

\begin{figure}[ptbh]
\centering
\begin{minipage}[t]{0.50\linewidth}
\centering
\includegraphics[scale=0.36]{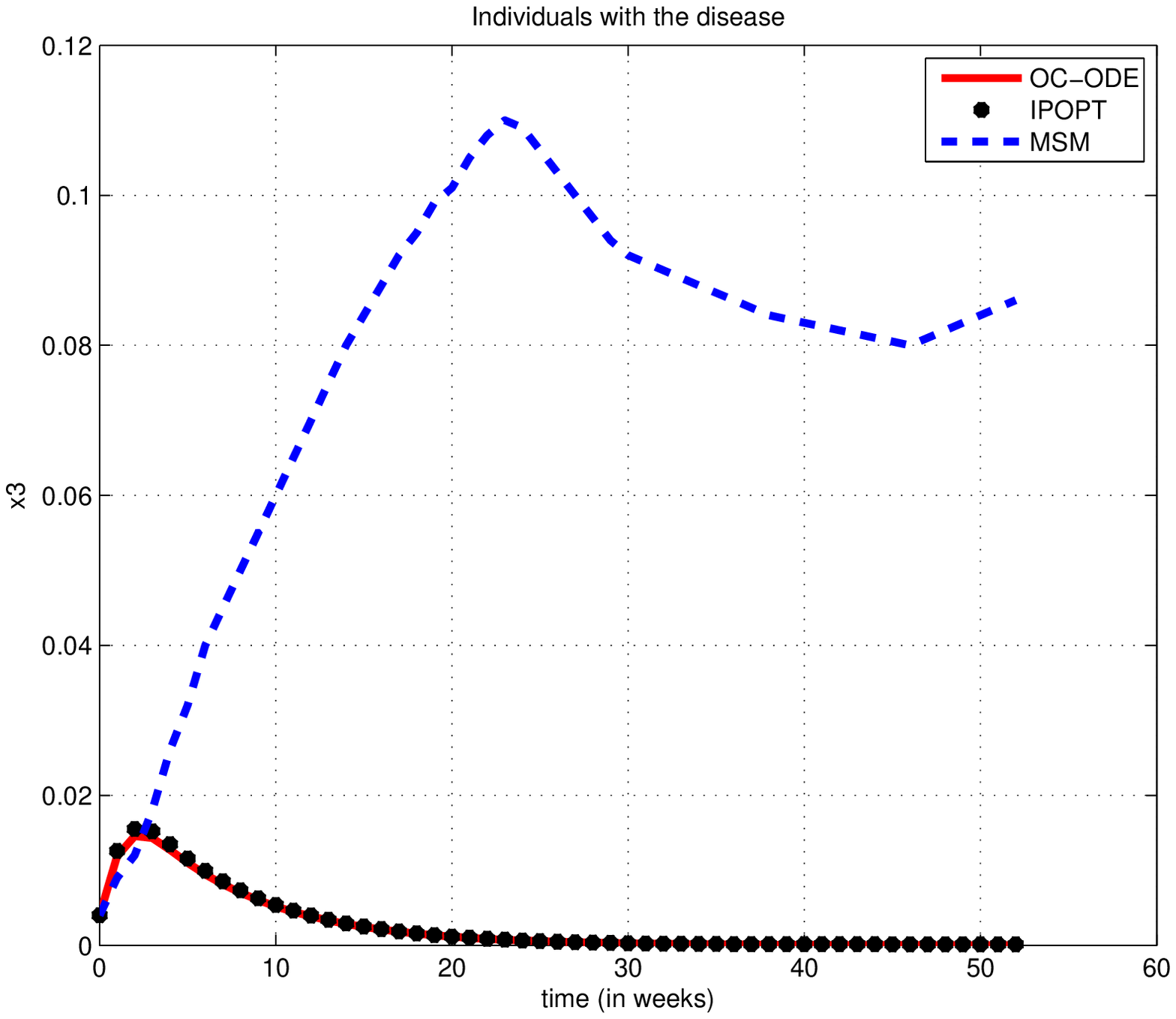}
{\caption{\label{plotx3} Infected individuals}}
\end{minipage}\hspace*{\fill}
\begin{minipage}[t]{0.50\linewidth}
\centering
\includegraphics[scale=0.36]{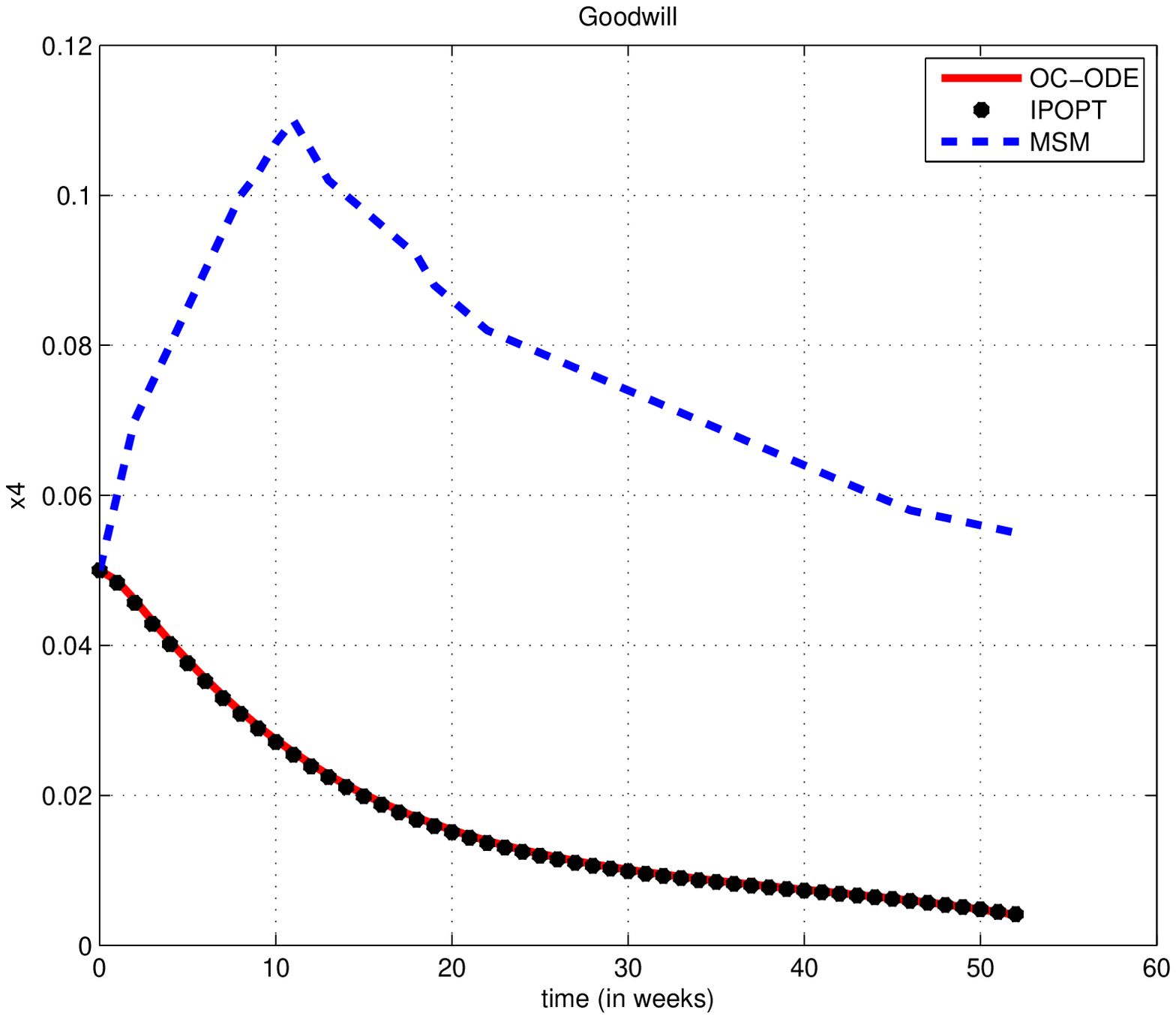}
{\caption{\label{plotx4} \small Level of popular motivation}} \label{pplot}
\end{minipage}
\end{figure}

\begin{figure}[ptbh]
\begin{center}
\includegraphics[scale=0.45]{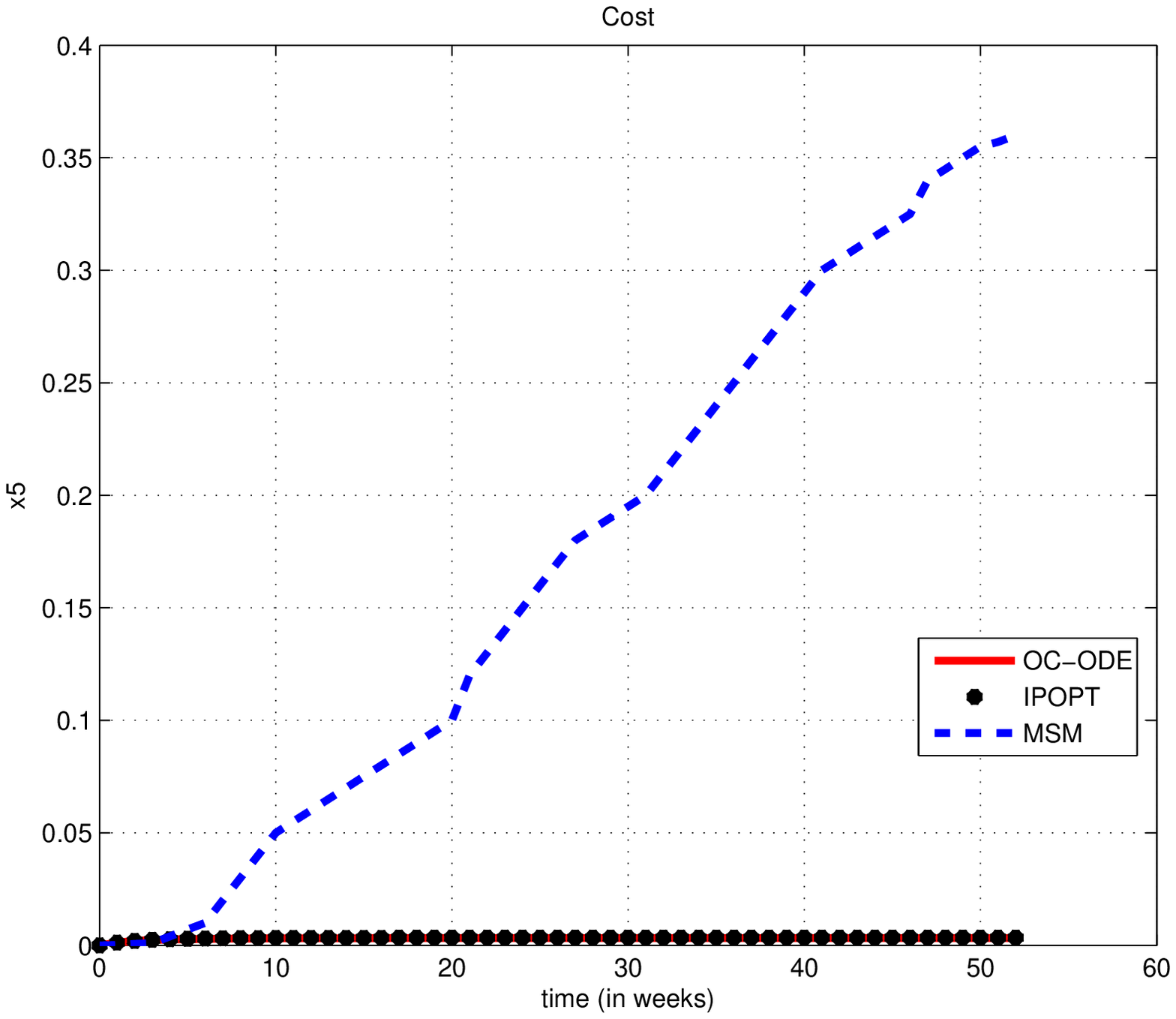}
\end{center}
\caption{Cost \label{plotx5}}
\end{figure}

\begin{figure}[ptbh]
\begin{minipage}[t]{0.50\linewidth}
\centering
\includegraphics[scale=0.36]{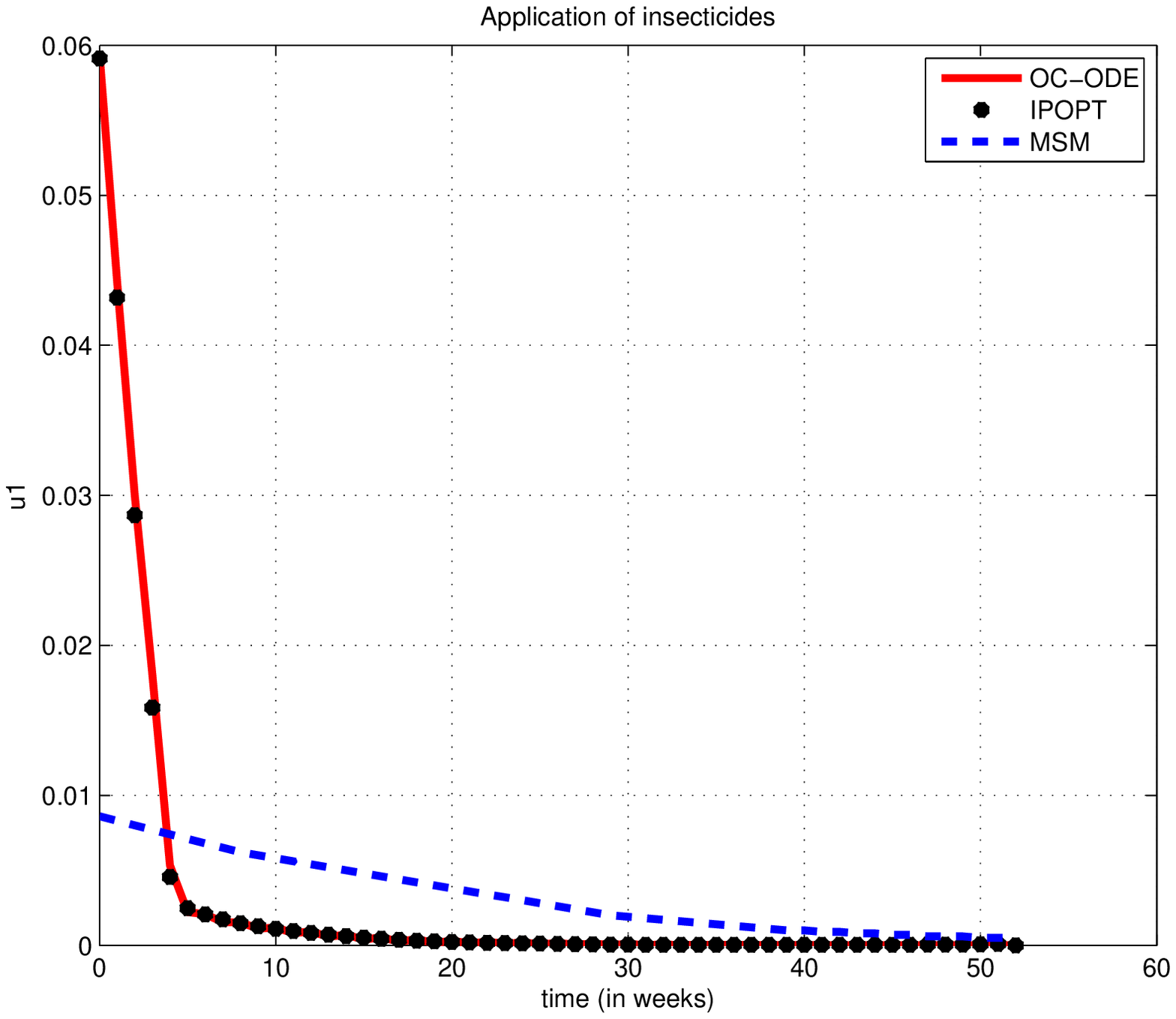}
{\caption{\label{plotu2} \small Application of insecticides}} \label{ppplot}
\end{minipage}\hspace*{\fill}
\centering
\begin{minipage}[t]{0.50\linewidth}
\centering
\includegraphics[scale=0.36]{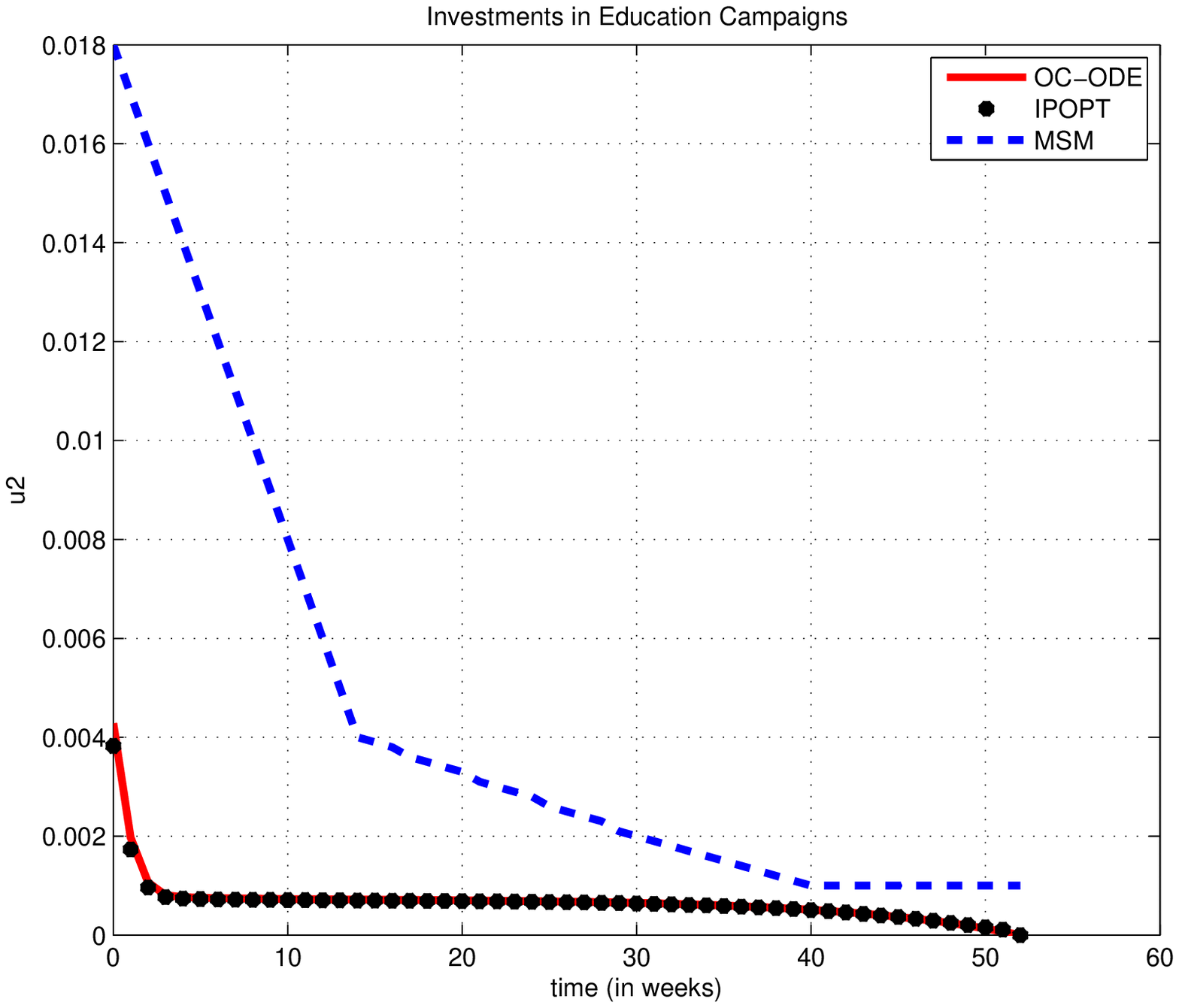}
{\caption{\label{plotu1}  Educational campaigns}}
\end{minipage}
\end{figure}

The results for the state and control variables are shown in Figures \ref{plotx1} to \ref{plotu2}.
Each figure has three graphics: \textsf{OC-ODE} and \textsf{IPOPT}, which correspond to the solutions obtained
by the solvers used, respectively; and MSM, corresponding to the Multiple Shooting Method \cite{Pesh,Trelat}
that was used by the authors of the paper \cite{Caetano2001}. It is important to salient that,
at the time of the initial paper \cite{Caetano2001}, the authors hadn't the same computational resources
that exist nowadays. The results we obtain using described methods
(\textsf{OC-ODE} and \textsf{IPOPT}) are better since the cost to combat
the Dengue disease and the number of infected individuals are smaller.

Figures~\ref{plotx1} and \ref{plotx2} show the density of Dengue mosquitoes.
It is possible to see that in this new solution, with the same number of mosquitoes
than in the previous solution \cite{Caetano2001},
the number of infected mosquitoes falls dramatically.
Figures~\ref{plotx3} and \ref{plotx4} report to the population in the risk area.
Our solution shows that the number of ill people decrease quickly.
That is also an explanation for the level of motivation to combat the mosquitos
to be also lower than the previous solution proposed in \cite{Caetano2001}.
Figure~\ref{plotx5} shows the accumulated cost. It is clear that almost all year
the cost is lower when compared with MSM \cite{Caetano2001}. This lower cost level
is a consequence of infected mosquitos and infected individuals both falling down under our approach.
Figures~\ref{plotu1} and \ref{plotu2} are related to the controls applied:
educational campaigns and application of insecticides. It is possible
to see that the new functions for the control variables are more economic.


\section{Conclusions}
\label{sec:conc}

At this moment the world, as a result of major demographic changes,
rapid urbanization on a massive scale, global travel and environmental change,
faces enormous future challenges from emerging infectious diseases.
Dengue illustrates well these challenges \cite{WHO}.

In this work we investigated an optimal control model to Dengue epidemics
proposed in \cite{Caetano2001}, that includes the dynamics of the Dengue mosquitoes,
the effect of educational campaigns, and where the cost functional reflects
a compromise between financial spending in insecticides and educational campaigns
and the population health. For comparison reasons, the same choice of data/parameters
as done in \cite{Caetano2001} was considered. However, the methods used in the paper
work also for other choices of the data and other problems \cite{Nosso1,Nosso2}.

Until some years ago, due to computational limitations, most of
the models were run using codes made by the authors themselves.
This was the case with \cite{Caetano2001}.
Nowadays, one can choose between several proper softwares ``out of the box'',
that already take into account specific features of stiff problems,
scaling problems, etc. With this work it is possible
to perceive that ``old'' problems can again be taken into account
and be better analyzed with new technology and approaches,
with the goal of finding global optimal solutions, instead of local ones.
The optimal control problem \cite{Caetano2001} was solved here numerically using
\textsf{OC-ODE} \cite{Matthias2009} and \textsf{IPOPT} \cite{Ipopt}.

The \textsf{OC-ODE} is a specific optimal control software,
while \textsf{IPOPT} is a standard nonlinear optimization solver.
It is possible to verify that, despite the different philosophies
of the \textsf{OC-ODE} and \textsf{IPOPT} solvers,
the reached is exactly the same. This fact enforces the robustness
of the obtained results. The first software uses
the optimal control theory in order to reach the solution.
The second one, uses nonlinear programming
as its main tool. It is important to mention
that the problem under study is a difficult one.
Other nonlinear packages were tested, and they could not reach
to a solution---some crashes at the middle or some bad scaling issues were observed.

The results obtained from \textsf{OC-ODE} and \textsf{IPOPT} coincide and improve the ones
previously reported in \cite{Caetano2001} (\textrm{cf.} Section~\ref{sec:4}).
Indeed, the control policy we obtain makes big progress with respect
to the previous best policy: the percentage of infected mosquitoes
vanishes just after four weeks, while mosquitoes are completely
irradicated after 30 weeks (Figures~\ref{plotx1} and \ref{plotx2});
the number of infected individuals begin to decrease after
four weeks while with the previous policy this happens after 23 weeks
only (Figure~\ref{plotx3}). Despite the fact that our results are better,
they are accomplished with a much smaller cost with insecticides
and educational campaigns (Figure~\ref{plotx5}).
The general improvement done, which explain why the results are so successful,
rely on an effective control policy of insecticides.
The proposed use of insecticides seems to explain
the big discrepancies between the results
here obtained and the best policy of \cite{Caetano2001}.
Our results show that applying insecticides in the first four weeks
yields a substantial reduction in the cost of combating Dengue,
in terms of the functional proposed in \cite{Caetano2001}.
The main conclusion is that health authorities should pay
attention to the epidemiology from the very beginning:
effective control decisions in the first four weeks have
a decisive role in the combat of Dengue, and
population and governments will both profit from it.

As future work we intend to analyze the model with different
parameters in the objective function and understand how the heights
associated to the variables can, or not, influence the decrease of
the disease. It will be also interesting to analyze other
problems, not only those related with the Dengue, such as
\cite{Duque2006}, but also with different biological illnesses,
such as \cite{Lenhart1997,ElGohary2009,Joshi2002}.


\section*{Acknowledgements}

The authors are grateful to the \emph{Portuguese Foundation for Science and Technology} (FCT)
for all the support: Rodrigues acknowledges the PhD grant SFRH/BD/33384/2008;
Monteiro the support of the \emph{Nonlinear Systems Optimization and Statistics group}
of the R\&D center ALGORITMI; Torres the support of the
\emph{Center for Research and Development in Mathematics and Applications} (CIDMA),
and the project UTAustin/MAT/0057/2008.
They would like to express their gratitude to Matthias Gerdts, for
having provided a copy of his software \cite{Matthias2009} as well as for
helpful discussions on the subject, and to two anonymous referees
for encouragement words and several relevant and stimulating remarks
contributing to improve the quality of the paper.




\begin{thebibliography}{99}

\bibitem{Betts}
    {\sc J. Betts},
    {\em Practical Methods for Optimal Control Using Nonlinear Programming},
    {\rm SIAM: Advances in Design and Control}, 2001.

\bibitem{Lenhart1997}
    {\sc S. Butler, D. Kirschner and S. Lenhart},
    {\em Optimal control of chemotherapy affecting the infectivity of HIV},
    Arino, O. (ed.) et al., Advances in mathematical population
    dynamics - molecules, cells and man. Papers from the 4th international
    conference, Rice Univ., Houston, TX, USA, May 23--27, 1995. Singapore: World
    Scientific Publishing. Ser. Math. Biol. Med. {\bf 6} (1997) 557--569.

\bibitem{Caetano2001}
    {\sc M. A. L. Caetano and T. Yoneyama},
    {\rm Optimal and sub-optimal control in Dengue epidemics},
    {\em Optimal Control Applications and Methods} {\bf 22} (2001) 63--73.

\bibitem{Duque2006}
    {\sc J. E. Duque and M. A. Navarro-Silva},
    {\rm Dynamics of the control of Aedes (Stegomyia) aegypti Linnaeus (Diptera, Culicidae)
    by Bacillus thuringiensis var israelensis, related with temperature,
    density and concentration of insecticide},
    {\em Rev. Bras. Entomol.} {\bf 50}-4 (2006) 528--533.

\bibitem{ElGohary2009}
    {\sc A. El-Gohary and I. A. Alwasel},
    {\rm The chaos and optimal control of cancer model with complete unknown parameters},
    {\em Chaos, Solitons \& Fractals} {\bf 42}-5 (2009) 2865--2874.

\bibitem{Edite}
    {\sc E. M. G. P. Fernandes},
    {\em Computa\c{c}\~{a}o Num\'{e}rica},
    {\rm Universidade do Minho, Braga}, 1998.

\bibitem{AMPL}
    {\sc  R. Fourer, D. M. Gay and B. W. Kernighan},
    {\em AMPL: A Modeling Language for Mathematical Programming},
    {\rm Duxbury Press / Brooks/Cole Publishing Company, 2002}.

\bibitem{Matthias2009}
    {\sc M. Gerdts},
    {\rm User's Guide OC-ODE (version 1.4)}
    {\em Technical Report}, Universit\"{a}t W\"{u}rzburg, May, 2009.

\bibitem{Joshi2002}
    {\sc H. R. Joshi},
    {\rm Optimal control of an HIV immunology model},
    {\em Optim. Contr. Appl. Math.} {\bf 23} (2002) 199--213.

\bibitem{Lewis1995}
    {\sc F. L. Lewis and V. L. Syrmos},
    {\em Optimal Control (2nd ed)},
    {\rm Wiley: New  York, 1995.}

\bibitem{NEOS}
{\sc NEOS}.
\newblock Webpage: \url{http://www-neos.mcs.anl.gov/neos/},
{\em October, 2009}.

\bibitem{Nerlove1962}
    {\sc M. Nerlove and K. J. Arrow},
    {\rm Optimal advertising policy under dynamic conditions},
    {\em Economica} {\bf 42}-114 (1962) 129--142.

\bibitem{Pesh}
    {\sc H. J. Pesch},
    {\rm Real time computation of feedback controls for constrained optimal control problems.
    Part 2: a correction method based on multiple shooting},
    {\em Optimal Control Applications and Methods} {\bf 10} (1989) 147--171.

\bibitem{Nosso1}
    {\sc H. S. Rodrigues, M. T. T. Monteiro and D. F. M. Torres},
    {\rm Optimization of Dengue epidemics: a test case with different discretization schemes},
    {\em American Institute of Physics Conf. Proc.} {\bf 1168} (2009) 1385--1388.
    {\tt arXiv:1001.3303}

\bibitem{Nosso2}
    {\sc H. S. Rodrigues, M. T. T. Monteiro, D. F. M. Torres and A. Zinober},
    {\rm Control of Dengue disease: a case study in Cape Verde},
    Proc. 10th International Conference on Mathematical Methods in Science and Engineering,
    Almer\'{\i}a 26-30 June 2010.

\bibitem{Trelat}
    {\sc E. Tr\'{e}lat},
    {\em Contr\^{o}le optimal: th\'{e}orie \& applications},
    {\rm Vuibert, Collection Math\'{e}matiques Concr\`{e}tes}, 2005.

\bibitem{Ipopt}
    {\sc A. W\"{a}chter and L. T. Biegler},
    {\rm On the implementation of an interior-point filter
    line-search algorithm for large-scale nonlinear programming},
    {\em Mathematical Programming: Series A and B.} {\bf 106}-1 (2006) 25--57.

\bibitem{WHO}
    \url{http://www.who.int/topics/dengue/en/}, October 2009.

\end{thebibliography}
\end{document}